\documentclass[10pt]{amsart}
\usepackage{latexsym}
\usepackage{amsmath}
\usepackage{amssymb}
\usepackage[all]{xy}

\newcommand{\C}{{\mathbb C}}
\newcommand{\Q}{{\mathbb Q}}
\newcommand{\Z}{{\mathbb Z}}
\newtheorem{thm}{Theorem }[section]
\newtheorem{prop}[thm]{Proposition}
\newtheorem{lem}[thm]{Lemma}

\newtheorem{cor}[thm]{Corollary}
\newtheorem{exmp}[thm]{Example}

\begin{document}
{

\title{A relaxed evaluation  subgroup}
\author{Toshihiro YAMAGUCHI}
\date{}
\address{Faculty of Education, Kochi University, 2-5-1, Kochi,780-8520, JAPAN}
\email{tyamag@kochi-u.ac.jp}
\maketitle

\begin{abstract}
Let $f:X\to Y$  be a pointed map between connected 
 CW-complexes.
As a generalization of the evaluation subgroup ${ G}_*(Y,X;f)$,
we will define the {\it relaxed evaluation  subgroup} 
${\mathcal G}_*(Y,X;f)$ in   the homotopy group 
$\pi_*(Y)$ of $Y$,
which is identified with ${\rm Im} \pi_*(\tilde{ev})$ for
the evaluation map
$\tilde{ev} :map(X,Y;f)\times X\to Y$ given by
$\tilde{ev} (h,x)=h(x)$.
Especially we see by using Sullivan model in rational homotopy theory
for the rationalized map $f_{\Q}$
  that ${\mathcal G}_*(Y_{\Q},X_{\Q};f_{\Q})=\pi_*(Y)\otimes \Q$ if 
the map $f$ induces an injection of rational homotopy groups.
Also we compare it with more relaxed subgroups
 by 
several  rationalized  examples.
\end{abstract}

\section{Introduction}

Let $f: X\to Y$ and $g: B\to Y$ be pointed maps of connected 
CW complexes.
Recall the definition of pairing with axes $f$ and $g$ \cite{O},
which is given by the existence of a map $F_{g,f}:B\times X\to Y$ 
in the homotopy commutative diagram:
$$(1)\ \ \ \ 
\xymatrix{B\times X
\ar@{.>}[dr]^{F_{g,f}}& X\ar[l]_{\ \ \ i_X}\ar[d]^{f}\\
B \ar[r]_{g}\ar[u]^{i_B}
 & Y, \\
}$$
where $i_B(b)=(b,*)$ and $i_X(x)=(*,x)$.
In particular, when $B=Y=X$ and $f=g=id$,
$X$ is an H-space of the multiplication $F_{id,id}$.
In this paper, we consider  whether or not 
 there exist 
 a section $s:B\to B\times X$ 
 and a map $F:B\times X\to Y$ such that
the  diagram 
$$(2)\ \ \ \ 
\xymatrix{B\times  X
\ar@{.>}[dr]^{F}& X\ar[l]_{\ \ i_X}\ar[d]^{f}\\
B \ar[r]_{g}\ar@{.>}[u]^s
 & Y. \\
}$$ 
homotopically commutes.
Here a section means a map $s:B\to B\times X$ which  satisfies 
$p_B\circ s\simeq id_B$ for the projection 
$p_B:B\times X\to B$ with $p_B(b,x)=b$.


The $n$th Gottlieb group 
$G_n(X)$ of a space $X$ is the subgroup of the 
$n$th homotopy group $\pi_n(X)$ of $X$
consisting of homotopy classes of maps $a:S^n\to X$ such that
the wedge $(a|id_X):S^n\vee X\to X$ extends
to a map $F_a:S^n\times X\to X$  \cite{G}.
The $n$th evaluation subgroup $G_n(Y,X;f)$ of a map $f:X\to Y$
is the subgroup of $\pi_n(Y)$ represented by  maps $a:S^n\to Y$ such that
 $(a|f):S^n\vee X\to Y$ extends
to a map $F_{a,f}:S^n\times X\to Y$  inducing the homotopy
commutative diagram  
$$(3)\ \ \ \ 
\xymatrix{S^n\times X
\ar@{.>}[dr]^{F_{a,f}}& X\ar[l]_{\ \ \ i_X}\ar[d]^{f}\\
S^n \ar[r]_{a}\ar[u]^{i_{S^n}}
 & Y \\
}$$
\cite{WK},
which is the case of $B=S^n$ in $(1)$.
The map $F_{a,f}$ in $(3)$ is 
the adjoint of 
a map $F_{a,f}'$ in the homotopy commutative diagram
$$(3)'\ \ \ \  
\xymatrix{
& map_f(X,Y)\ar[d]^{ev}\\
S^n \ar@{.>}[ur]^{F_{a,f}'}\ar[r]_{a}
 & Y \\
}$$
where
$F'_{a,f}(b)(x):=F_{a,f}(b,x)$ for $b\in S^n$ and $x\in X$.
Here $map_f(X,Y)$ is the connected component of $f$ in
the mapping space $map(X,Y)$ from $X$ to $Y$
with the compact open topology
and $ev$ is the evaluation map
given by $ev(h)=h(*)$. 
Then $$G_n(Y,X;f)={\rm Im} \ (\pi_n(ev):\pi_n(map_f(X,Y))\to \pi_n(Y))$$ in $\pi_n(Y)$
and therefore it is called an `evaluation' subgroup of a map.
Notice that 
 $(3)$ is a special case of the homotopy commutative diagram
 $$(4)\ \ \ \  
\xymatrix{S^n\times X
\ar@{.>}[dr]^{F}& X\ar[l]_{\ \ i_X}\ar[d]^{f}\\
S^n \ar[r]_{a}\ar@{.>}[u]^s
 & Y,\\
}$$
in which 
$\eta:X\overset{i_X}{\to} S^n\times X\overset{p_{S^n}}{\to} S^n$
is  a trivial  fibration
with a section $s:S^n\to S^n\times X$. 
Here we can put $s(b)=(b,\tau (b))$ for a map $\tau :S^n\to X$.\\

\noindent
{\bf Definition A.} 
For a  map $f:X\to Y$, the $n$th {\it relaxed evaluation 
 subgroup} of $f$  is given as 
$${\mathcal G}_n(Y,X;f):=\{ a\in \pi_n(Y)\ |\  there \ are \  a \ section \ s:S^n\to S^n\times X$$ 
$$and 
\ a \ 
map \ F:S^n\times X\to Y\
 such \ that\ 
\ F\circ s\simeq a,\  f\simeq F\circ i_X \}.$$ 
\vspace{0.2cm}

The map $F$ in $(4)$ 
is the adjoint of a map $\tilde{F'}$ in the homotopy commutative diagram
$$(4)'\ \ \ \  
\xymatrix{
& map_f(X,Y)\times X\ar[d]^{\tilde{ev}}\\
S^n \ar@{.>}[ur]^{\tilde{F'}}\ar[r]_{a}
 & Y \\
}$$
where 
$\tilde{F'}(b):=(F'(b),\tau (b))$
with $F'(b)(x)=F(b,x)$ and $s(b)=(b,\tau (b))$ for $b\in S^n$ and $x\in X$.
Here $\tilde{ev}:map(X,Y;f)\times X\to Y$ is the evaluation map
given by $\tilde{ev}(h,x)=h(x)$.
Thus\\

\noindent
{\bf Claim.}  \ 
  ${\mathcal G}_n(Y,X;f)={\rm Im} (\pi_n(\tilde{ev}):\pi_n(map_f(X,Y)\times X)\to \pi_n(Y))$.\\

In  this paper,  we will estimate
${\rm Im} \pi_n(\tilde{ev})$
in several cases according to Definition A.
We note that
it is a naturally generalized object of an ordinary evaluation subgroup.
Indeed, for a subcomplex $X'$ of $X$,
we can put
$${\mathcal G}_n(Y,X;f)(X'):=\{ a\in \pi_n(Y)\ |\  there \ are \  a \ section \ s:S^n\to S^n\times X'$$ 
$$and 
\ a \ 
map \ F:S^n\times X\to Y\
 such \ that\ 
\ F\circ s\simeq a,\  f\simeq F\circ i_X \}.$$
Then we have $${\mathcal G}_n(Y,X;f)\supset {\mathcal G}_n(Y,X;f)(X')\supset {\mathcal G}_n(Y,X;f)(*)=
G_n(Y,X;f)$$
and ${\mathcal G}_n(Y,X;f)(X')={\rm Im} (\pi_n(\tilde{ev}):\pi_n(map_f(X,Y)\times X')\to \pi_n(Y))$.
In the following, we see several properties of a relaxed evaluation subgroup.

\begin{lem}
For 
a subspace $X$ of a space $Y$
and a map $g:B\to Y$ with ${\rm Im} \ g\subset X$,
there are a section $s:B\to B\times X$ and a map $F:B\times X\to X$
such  that the diagram 
$$
\xymatrix{B\times X
\ar@{.>}[dr]^{F}& X\ar[l]_{\ \ i_X}\ar[d]^{\cap}\\
B \ar[r]_{g}\ar@{.>}[u]^s
 & Y\\
}$$
 commutes. 
\end{lem}

Indeed, put  $s(b)=(b,g(b))$ and
$F((b,x))=x$
for $(b,x)\in B\times X$.

In particular,  in the case of $B=S^n$ 
and $X=Y^m$, the $m$-skelton of $Y$, we have
by  cellular approximation theorem

\begin{prop}\label{p}   
For $n\leq m$,
${\mathcal G}_n(Y,Y^m;i_{Y^m})= \pi_n(Y)$
for $i_{Y^m} : Y^m\subset Y$.
In particular,
${\mathcal G}_*(Y,Y;id_Y)= \pi_*(Y)$.
\end{prop}

Of course ${\mathcal G}_n(Y,*;*)=  G_n(Y,*;*)=\pi_n(Y)$
for the constant map $*:*\to Y$
and we know
 $G_n(Y)=G_n(Y,Y;id_Y)\subset G_n(Y,X;f)$ for any map $f:X\to Y$.
In contrast,
${\mathcal G}_n(Y,X;f)\subset {\mathcal G}_n(Y,Y;id_Y)$
from Proposition \ref{p}.
Note that ${\mathcal G}_*(Y,X;f)$ may be zero for a map  $f:X\to Y$
even if $\pi_*(Y)\neq 0$ (see Example \ref{c} and examples in Section 3).

Also  as an evaluation subgroup satisfies it, 
a relaxed evaluation subgroup satisfies
$$\pi_n(g)({\mathcal G}_n(Y,X;f))\subset  
{\mathcal G}_n(Y',X;g\circ f)$$
\noindent
for a pointed map $g:Y\to Y'$.
Thus there is a  map
$\pi_n(f):\pi_n(X)={\mathcal G}_n(X,X;id_X)\to 
{\mathcal G}_n(Y,X;f)$.
Therefore the relaxed eveluation  subgroups
${\mathcal G}_n(X,Z;j)$ and
${\mathcal G}_n(Y,X;f)$
are embedded into
parts of the homotopy exact sequence of a fibration
$\xi :Z\overset{j}{\to} X\overset{f}{\to} Y$, respectively. 
\begin{cor}
For a fibration
$\xi :Z\overset{j}{\to} X\overset{f}{\to} Y$ and any $n$,
there are the sequences 
$$
\pi_{n+1}(Y)\overset{\partial_{n+1}}{\to} 
\pi_n(Z)\overset{\pi_n(j)}{\to} 
{\mathcal G}_n(X,Z;j)\overset{\pi_n(f)}{\to}
 \pi_{n}(Y)$$
and 
$$ 
\pi_{n}(Z)\overset{\pi_n(j)}{\to} 
\pi_n(X)\overset{\pi_n(f)}{\to} 
{\mathcal G}_n(Y,X;f)\overset{\partial_n}{\to}
 \pi_{n-1}(Z)$$
which are both exact.
Moreover,
for the connecting map $\tilde{\partial} :\Omega Y\to Z$ of $\xi$,
$$ \pi_{n+1}(X)\overset{\pi_{n+1}(f)}{\to}
\pi_{n+1}(Y)\overset{\partial_{n+1}}{\to} 
{\mathcal G}_n(Z,\Omega Y;\tilde{\partial})\overset{\pi_n(j)}{\to}
 \pi_{n}(X)$$
is exact.
\end{cor}

 Note that, for a pointed map $g:X'\to X$,
there is an inclusion ${ G}_n(Y,X;f)\subset  
{ G}_n(Y,X';f\circ g)$.
But it does not hold for  relaxed evaluation subgroups
in general.

\begin{lem}\label{ss}
For a map $g:B\to Y$ and 
a map $f:X\to Y$ such that  $f_*:[B,X]\to [B,Y]$
is surjective,
there are a section $s:B\to B\times X$ and a map $F:B\times X\to X$
such  that the diagram 
$$
\xymatrix{B\times X
\ar@{.>}[dr]^{F}& X\ar[l]_{ \ \ i_X}\ar[d]^{f}\\
B \ar[r]_{g}\ar@{.>}[u]^s
 & Y\\
}$$
homotopically  commutes. 
\end{lem}

Indeed, there is a  lift $\tilde{g} :B\to X$
such that $f\circ \tilde{g}\simeq g$ from the assumption.
Then put  $s(b)=(b,\tilde{g} (b))$ and
$F((b,x))=f(x)$
for $(b,x)\in B\times X$.

In particular,  in the case of $B=S^n$, we have

\begin{prop}\label{s}  
If a map $f:X\to Y$ induces a surjection $\pi_n(f):\pi_n(X)\to
\pi_n(Y)$,
 then
${\mathcal G}_n(Y,X;f)= \pi_n(Y)$.
\end{prop}

For example, ${\mathcal G}_*(Y,X;f)= \pi_*(Y)$
if a map $f$ has a section.
In contrast, $G_*(Y,X;f)=G_*(Y)$ if $f$ has a section.

Suppose that 
$f_{\Q}:X_{\Q}\to Y_{\Q}$ is a rationalized map
with $X$ and $Y$  
simply connected CW complexes of finite type \cite{HMR}.
We consider the relaxed evaluation  subgroup of a map
and more relaxed subgroups from a point of view of rational homotopy.
By using Sullivan's minimal model arguments, we show

\begin{thm}\label{b}
If a map $f:X\to Y$ induces an injection $\pi_*(f)\otimes \Q:\pi_*(X)\otimes \Q\to
\pi_*(Y)\otimes \Q$ on rational homotopy groups,
 then ${\mathcal G}_*(Y_{\Q},X_{\Q};f_{\Q})=\pi_*(Y)\otimes \Q$.
\end{thm}

From Proposition \ref{s} and Theorem \ref{b}, we have

\begin{cor}
If the homotopy fiber  of a map $f:X\to Y$ has the rational homotopy type of 
the Eilenberg-MacLane space $K(\Q , n)$ for some $n$, then 
${\mathcal G}_*(Y_{\Q},X_{\Q};f_{\Q})=\pi_*(Y)\otimes \Q$.
\end{cor}

In Section 2, we prove Theorem \ref{b}
after preparing of the notion of derivations of model.
In Section 3, we will define more relaxed subgroups of the homotopy group of $Y$,
 the {\it tncz subgroup} ${\mathcal T}_*(Y,X;f)$ and 
the {\it sectional subgroup}  ${\mathcal S}_*(Y,X;f)$, for a map $f:X\to Y$.
They are defined by relaxing the trivial fibration $\eta$ in $(4)$.
 We compare ${\mathcal G}_*(Y,X;f)$ with them  by several rationalized  examples.



\section{Sullivan  models}

We use the {Sullivan minimal model} $M(Y)$ 
 of a simply connected CW complex $Y$ of finite type.
It is a free $\Q$-commutative differential graded algebra  (DGA) 
 $(\Lambda{W},d_Y)$
 with a $\Q$-graded vector space $W=\bigoplus_{i\geq 2}W^i$
 where $\dim W^i<\infty$ and a decomposable differential.
 Here  $\Lambda^+{W}$ is 
 the ideal of $\Lambda{W}$ generated by elements of positive degree.
Denote the degree of a homogeneous element $x$ of a graded algebra as $|{x}|$.
Then  $xy=(-1)^{|{x}||{y}|}yx$ and $d(xy)=d(x)y+(-1)^{|{x}|}xd(y)$. 
 Recall $M(Y)$ determines the rational homotopy type of $Y$.
Especially there is an isomorphism $W^i\cong Hom(\pi_i(Y),{\Q})$.
We denote the dual element of $a\in \pi_i(Y)\otimes \Q$ as $a^*$.
 Put $M(Y)=(\Lambda W,d_Y)$. 
  Then the model   of a map $f: X{\to} Y$ is given by 
a KS-extension $$(\Lambda W,d_Y)\overset{i}{\to}
 (\Lambda W\otimes \Lambda V,D)\overset{q}{\to} (\Lambda V,\overline{D})$$
 with  $D|_{\Lambda W}=d_Y$ and the minimal model 
 $ (\Lambda V,\overline{D})$ of the homotopy fiber of $f$
or $$(H^*(Y;\Q) ,0)\overset{i}{\to}
 (\Lambda W\otimes \Lambda V,D)\overset{q}{\to} (\Lambda V,\overline{D})$$
when $Y$ is formal (for example, when $Y$ is a sphere) \cite{FHT}. 
 In general, $D$ is not decomposable
and it is decomposable 
if and only if 
 $\pi_*(f)\otimes \Q:\pi_*(X)\otimes \Q\to
\pi_*(Y)\otimes \Q$ is a surjection.
See \cite{FHT} for a general introduction and the standard notations.

Let $A$ be a DGA $A=(A^*,d_A)$ with $A^*=\oplus_{i\geq 0}A^i,\ A^0=\Q$, $A^1=0$  and 
the augmentation $\epsilon:A\to \Q$.
Define  
 $Der_i A$  the vector space of self-derivations of $A$
decreasing the degree by $i>0$,
where  $\theta(xy)=\theta(x)y+(-1)^{i|x|}x\theta(y)$
for $\theta\in Der_iA$. 
We denote  $\oplus_{i>0} Der_iA$ by
$DerA$.
The boundary operator $\delta : Der_* A\to Der_{*-1} A$
is 
defined by $\delta (\sigma)=d_A\circ \sigma-(-1)^{|\sigma |}\sigma\circ d_A$.
For a DGA-map $\phi:A\to B$,
define a $\phi$-derivation of degree $n$ to be a linear map
$\theta :A^*\to B^{*-n}$
 with $\theta(xy)=\theta(x)\phi (y)+(-1)^{n|x|}\phi (x)\theta(y)$
 and 
$Der(A,B;\phi)$ the vector space of $\phi$-derivations. 
The boundary operator $\delta_{\phi}: Der_* (A,B;\phi )\to Der_{*-1} (A,B;\phi )$ is defined by $\delta_{\phi}
 (\sigma)=d_B\circ\sigma-(-1)^{|\sigma |}\sigma\circ d_A$. 
Note $Der_* (A,A;id_A)=Der_*(A)$.
For 
$\phi:A\to B$, 
the composition with $\epsilon':B\to \Q$ 
induces a chain map $\epsilon'_*:
 Der_n(A,B;\phi )\to Der_n(A,\Q;\epsilon )$.
 For a minimal model $A=(\Lambda Z,d_A)$, define 
$G_n(A,B;\phi):=Im (H(\epsilon'_*):H_n(Der(A,B;\phi))\to Hom_n(Z,\Q)).$
%
Especially  $G_*(A,A;id_A)=G_*(A)$. 
Note that $z^*\in Hom (Z,\Q)$ ($z^*$ is the dual  of the basis element $z$) is in 
 $G_n(A,B;\phi )$
  if and only if 
  $z^*$ extends to a derivation $\theta\in Der(A,B;\phi )$
  with $\delta_{\phi} (\theta)=0$.
\begin{thm}{\rm(\cite{FH},\cite{LS})}\ \ 
$ G_n(Y_{\Q},X_{\Q};f_{\Q})\cong G_n(M(Y),M(X);M(f))=
G_n((\Lambda W,d_Y), (\Lambda W\otimes \Lambda V,D))$.
\end{thm}
Thus $G_n(Y_{\Q},X_{\Q};f_{\Q})$ is completely determined
by the derivations only. 
For a DGA-map $\phi :(\Lambda V,d)\to (\Lambda Z,d')$, 
 the symbol $(v,h)\in Der (\Lambda V,\Lambda Z;\phi )$
 means the $\phi$-derivation sending an element $v\in V$
to $h\in \Lambda Z$ and the other to zero.
Especially $(v,1)=v^*$.

\begin{exmp}\label{nos}{\rm 
Consider the fibration $S^3\to X\overset{f}{\to} Y=S^2\times S^2$
whose KS-extension is given by
$$(\Lambda (w_1,w_2,w_3,w_4),d_Y)\to (\Lambda (w_1,w_2,w_3,w_4,v),D)
\to (\Lambda (v),0),$$
where $|w_1|=|w_2|=2$, $|w_3|=|w_4|=|v|=3$,
$d_Yw_1=d_Yw_2=0$, $d_Yw_3=w_1^2$, $d_Yw_4=w_2^2$
and $Dv=w_1w_2$.
Since $D$ is decomposable,  
$\pi_*(f)\otimes \Q:\pi_*(X)\otimes \Q\to \pi_*(Y)\otimes \Q$
is surjective.
So we have $
{\mathcal G}_2(Y_{\Q},X_{\Q};f_{\Q})=\pi_2(Y)\otimes \Q$ from 
Proposition  \ref{s}.
On the other hand, 
from Theorem 2.1, $
{ G}_2(Y_{\Q},X_{\Q};f_{\Q})=0$
since $$\delta_f((w_i,1))=2(w_{i+2},w_i)\not\in \delta_f\ (Der 
(\Lambda W,\Lambda W\otimes \Lambda^+ v))$$
for $i=1,2$ and $W=\Q (w_1,w_2,w_3,w_4)$.
Note that $f_{\Q}$ has no section (\cite{T}).
}
\end{exmp}

\vspace{0.2cm}

\noindent
{\it Proof of Theorem \ref{b}}. 
 Fix an element $a\in \pi_n(Y)\otimes \Q$.
From the assumption,
there is a DGA-projection
$p_{W,V}:(\Lambda W,d_Y)\to (\Lambda V,\overline{d}_Y)$
as the model of $f$.
 A model of the (non homotopy commutative) diagram
 $$ \ 
\xymatrix{
E\ar[d]_{p}
 & 
X_{\Q}\ar[l]_{i}\ar[d]^{f_{\Q}}  \\
S^n_{\Q}\ar[r]_{a}  &  
Y_{\Q}\\
}$$
 is given by 
$$
\xymatrix{
(\Lambda x/x^2\otimes  \Lambda V,{D'}) 
\ar[r]^{\ \ \ \ \ q} & 
( \Lambda V,\overline{d}_Y)  \\
(\Lambda x/x^2,0) \ar[u]^{\cup} &  
( \Lambda W,d_Y)\ar[u]_{p_{W,V}}\ar[l]^{M(a)}\\
}$$
where $|x|=n$, $i_{V,W}:V\subset W$
and $p_{W,V}(\Lambda W)=\overline{\Lambda W}=\Lambda V$
with $p_{W,V}\circ i_{V,W}=id_V$.
Here $\Lambda x/x^2=\Lambda x$ if $n$ is odd and
 $\Lambda x/x^2=\Q[ x]/(x^2)$ if $n$ is even.

 We will construct a rationally trivial  fibration
 of the form $\eta :X_{\Q}\to E\to S^n_{\Q}$
 together with a suitable map $F:E\to Y_{\Q}$,
 in  model terms.
Define a graded $\Q$-algebra map 
$F':\Lambda W\to \Lambda x/x^2\otimes \Lambda V$ by
$$F'(w)=\overline{w}+(-1)^{n|w|}\sigma (w)x$$
where $\sigma \in Der_n(\Lambda W,  \Lambda V ;p_{W,V})$ 
 with $(-1)^{n|u|}\sigma (u) x=M(a)(u)$
for $ u\in W$. 
Also define the differential $D'$  by $D'(x)=0$ and 
$$D'|_{\Lambda V}=\overline{d}_Y-\delta_{\overline{d}_Y}(\overline{\sigma} ) \cdot x,$$
where 
$\overline{\sigma}
\in Der_{n}(\Lambda V)$
is uniquely given  by 
$\overline{\sigma}(\overline{w})=\sigma (w)$ for $w\in \Lambda W$
and $(\theta\cdot x)(z):=(-1)^{n|z|}\theta (z)x$ for a derivation $\theta$.
Then $D'\circ D'=0$ from $\delta_{\overline{d}_Y}\circ \delta_{\overline{d}_Y}=0$ 
and $F'$ is a DGA-map by
$$F'd_Y(w)=\overline{d_Y(w)}
+(-1)^{n(|w|+1)}\sigma(d_Yw)x$$
$$=\overline{d_Y(w)}
-(-1)^{n|w|}\delta_{\overline{d}_Y}(\overline{\sigma} )(\overline{w})x
+(-1)^{n|w|}\overline{d}_Y\sigma(w) x$$
$$=D'(\overline{w}+(-1)^{n|w|}\sigma (w)x)=D'F'(w)$$
for $w\in \Lambda W$.
Thus we have the KS-model of $\eta$
 $$(\Lambda x/x^2,0)\overset{i}{\to}
 (\Lambda x/x^2\otimes \Lambda V,D')\overset{q}{\to} ( \Lambda V,{\overline{d}_Y})$$
 and a map
 $$F': ( \Lambda W,d_Y)\to (\Lambda x/x^2\otimes  \Lambda V,{D'}). $$
Since $ \delta_{\overline{d}_Y}(\overline{\sigma} )\in Der(\Lambda V,
\Lambda^+   V)$,
the fibration $\eta $ 
has a section $s$
(\cite{T}).
Moreover the  definition of $D'$ indicates
the rational triviality
of $\eta$: 
$$  (\Lambda x/x^2\otimes  \Lambda V,{D'})\cong (\Lambda x/x^2,0)
\otimes ( \Lambda V,{\overline{d}_Y})
$$
over $ (\Lambda x/x^2,0)$
since then the homotopy class of the classifying map
$S^n\to Baut_1X_{\Q}$, 
$[\delta_{\overline{d}_Y}(\overline{\sigma} )]$, is zero in 
 $\pi_n(Baut_1X)\otimes \Q=H_{n-1}(Der M(X))$ \cite{S}.
We can choose the model of $s$ 
by $M(s)(x)=x$ and $M(s)(z)=0$ for $z\in  \Lambda V$,
then $M(s)\circ F'= M(a)$.
Thus there is a DGA-commutative  diagram  
$$
\xymatrix{
(\Lambda x/x^2\otimes\Lambda V,{D'}) 
\ar[r]^{\ \ \ \ \ q}\ar[d]_{M(s)} & 
( \Lambda V,{\overline{d}_Y})  \\
(\Lambda x/x^2,0)  &  
( \Lambda W,d_Y).\ar[u]_{p_{W,V}}\ar[l]^{M(a)}\ar@{.>}[ul]^{F'}\\
}$$
It is the rational model of $(4)$.
\hfill\qed\\

\noindent
{\it Proof of Corollary 1.7}. 
Put the model of homotopy fiber $(\Lambda v,0)$.
When $Dv$ is decomposable,
we have from Proposition \ref{s}.
Also when $Dv$ is not decomposable,
there is a surjection $M(Y)=(\Lambda W,d_Y)\to M(X)=(\Lambda V,d_X)$
with
$\dim V=\dim W-1$.
Then we have from Theorem \ref{b}.
\hfill\qed\\

\noindent
{\bf Remark 1.}
 In the proof of above,
 the fibration $\eta :X_{\Q}\overset{i}{\to} E\overset{p}{\to} S^n_{\Q}$ is trivial, that is,
 there is a homotopy commutative diagram
 $$ \  
\xymatrix{E\ar[d]_p
& X_{\Q}\ar[l]_{i}\ar[d]^{i_{X_{\Q}}}\\
S^n_{\Q} 
 & S^n_{\Q}\times X_{\Q}\ar@{.>}[ul]_{g}\ar[l]^{p_{S^n_{\Q}}} \\
}$$
where 
$g$ is a homotopy equivalence.
The model
$M(g)$ is given by $id-\overline{\sigma}\cdot x$,
which
 is a quasi-isomorphism.
But  the map $g$ does not induce 
 the  homotopy commutative diagram
 with a section $s$ of $p$
 $$ \  
\xymatrix{E
& X_{\Q}\ar[l]_{i}\ar[d]^{i_{X_{\Q}}}\\
S^n_{\Q} \ar[r]_{i_{S^n_{\Q}}}\ar[u]^s
 & S^n_{\Q}\times X_{\Q}\ar@{.>}[ul]_{g} \\
}$$
in general.
Therefore, even if $\pi_*(f)\otimes \Q$ is surjective or injective,
we can not induce 
${G}_*(Y_{\Q},X_{\Q};f_{\Q})= \pi_*(Y)\otimes \Q$
in general.
For example, see Example \ref{nos} or Example \ref{d},  respectively.\\

In the followings,
we consider some examples, 
 whose models are given by  $M(s)(x)=x$ and
$M(s)(z)=0$ for $z\in M(X)$ as
in the proof of Theorem \ref{b}.   
The index of an element means the degree.

\begin{exmp}{\rm For $n>0$,
${\mathcal G}_n(S^n_{\Q},S^n_{\Q};id_{S^n_{\Q}})=\pi_n(S^n)\otimes \Q$ 
is given by
the following  commutative  diagrams 
with $|x|=n$.
$$
{\rm ({\it n}:odd)}\ \ \ 
\xymatrix{
(\Lambda x/x^2\otimes  \Lambda  w_n
,{0}) 
\ar[r]^{\ \ \ \ \ q}\ar[d]_{M(s)} & 
(\Lambda w_n,0)  \\
(\Lambda x/x^2,0)  &  
( \Lambda w_n,0)\ar[u]_{=}\ar[l]^{M(a)}\ar@{.>}[ul]^{F}\\
}$$
where 
$M(s)(x)=x$, $M(s)(w_n)=0$  and
$F(w_n)=w_n+cx$ for $M(a)(w_n)=cx$ with $c\in \Q$.
$$
{\rm ({\it n}:even)}\ \ \ 
\xymatrix{
(\Lambda x/x^2\otimes  \Lambda ( w_n,w_{2n-1})
,{D'}) 
\ar[r]^{\ \ \ \ \ q}\ar[d]_{M(s)} & 
(\Lambda (w_n,w_{2n-1}),d_Y)  \\
(\Lambda x/x^2,0)  &  
( \Lambda (w_n,w_{2n-1}),d_Y)\ar[u]_{=}\ar[l]^{M(a)}\ar@{.>}[ul]^{F}\\
}$$
where $d_Yw_n=0$, $d_Yw_{2n-1}=w_n^2$, $D'w_n=0$,  $D'w_{2n-1}=w_n^2+2cw_nx$,
$F(w_n)=w_n+cx$ and $F(w_{2n-1})=w_{2n-1}$.

}
\end{exmp}

\begin{exmp}{\rm 
For the Hopf fibration $S^3\to S^7\overset{f}{\to} S^4$,
we know ${G}_4(S^4_{\Q},S^7_{\Q};f_{\Q})=\pi_4(S^4)\otimes \Q$ \cite{LS2}.
Indeed, the KS-extension is given by
$$(\Lambda (w_4,w_7),d_Y)\to (\Lambda (w_4,w_7,v_3),D)\to (\Lambda (v_3),0)$$
with $d_Yw_4=0$, $d_Yw_7=w_4^2$ and $Dv_3=w_4$.
Then, from Theorem 2.1, 
${G}_4(S^4_{\Q},S^7_{\Q};f_{\Q})=G_4((\Lambda (w_4,w_7),d_Y),(\Lambda (w_4,w_7,v_3),D))=
\Q(w_4^*)$ since
 $\delta_f((w_4,1)-2(w_7,v_3))=0$. 
Thus we have
$${G}_4(S^4_{\Q},S^7_{\Q};f_{\Q})={\mathcal G}_4(S^4_{\Q},S^7_{\Q};f_{\Q})=\pi_4(S^4)\otimes \Q\cong \Q.$$

Also for the product 
fibration $S^3\to  X=S^7\times S^4\overset{f\times id}{\to} S^4\times S^4=Y$
with the trivial fibration
$*\to S^4\to S^4$,
we have  ${G}_4(Y_{\Q},X_{\Q};(f\times id)_{\Q})\cong \Q$ and $
{\mathcal G}_4(Y_{\Q},X_{\Q};(f\times id)_{\Q})=\pi_4(Y)\otimes \Q\cong \Q\oplus \Q$.
}
\end{exmp}

\begin{exmp}\label{d}{\rm
Suppose that  the model of a map $f:X\to Y$ is given by
the projection $$p_{W,V}:M(Y)=(\Lambda (w_3,w_5,w_7,w_9),d_Y)\to  (\Lambda (w_3,w_7,w_9),\overline{d}_Y)=M(X)$$
where $p_{W,V}(w_5)=\overline{w_5}=0$,
$d_Yw_3=d_Yw_5=0$, $d_Yw_7=w_3w_5$, $d_Yw_9=w_3w_7$
and $\overline{d}_Yw_3=\overline{d}_Yw_7=0$
and $\overline{d}_Yw_9=w_3w_7$.
Then we have 
$G_3(Y_{\Q},X_{\Q};f_{\Q})
=G_7(Y_{\Q},X_{\Q};f_{\Q})=0$ by the direct calculations of derivations. 
But we have
${\mathcal G}_*(Y_{\Q},X_{\Q};f_{\Q})=\pi_*(Y)\otimes \Q$
from Theorem \ref{b}.
}
\end{exmp}

\begin{exmp}\label{c}
{\rm When $f:S^3\times S^3\to S^6$ is   the map  collapsing $S^3\vee S^3$, 
then we show
$${\mathcal G}_6(S^6,S^3\times S^3;f)=0.$$

 By degree arguments,
 any fibration $\eta :S^3\times S^3
\overset{i}{\to} E\overset{p}{\to} S^6$  is 
 rationally trivial, especially  $D=0$ in the KS-extension
$$(\Q [w_6]/(w_6^2),0)\to  (\Q [w_6]/(w_6^2)\otimes \Lambda (u_3,v_3),D)\to (\Lambda (u_3,v_3),0),$$
where $H^*(S^6;\Q)=\Q [w_6]/(w_6^2)$
and $M(S^3\times S^3)=(\Lambda (u_3,v_3),0)$.
In particular  $E_{\Q}\simeq 
(S^6\times S^3\times S^3)_{\Q}$.

For  $a\neq 0\in \pi_6(S^6)$,
suppose that there exists a map $F:E\to S^6$ 
in
 $$
\xymatrix{E
\ar@{.>}[dr]^{F}& S^3\times S^3\ar[l]_{i}\ar[d]^{f}\\
S^6 \ar[r]_{a}\ar@{.>}[u]^s
 & S^6,\\
}$$

Then
 $$M(F_{f,g})(w_6)=cx+u_3v_3$$
 for some non-zero $c\in \Q$
 associated to $a$.
 Then $w_6^2=0$ in $H^{12}(S^6;\Q)$ but
 $$[M(F_{f,g})(w_6)]^2=[2cxu_3v_3]\neq 0$$
 in $H^{12}(E;\Q)=H^{12}(S^6\times S^3\times S^3;\Q)$.
 It is a contradiction.
}
\end{exmp}

\section{The other subgroups}

Futher we will relax Definition A.
For pointed maps $f:X\to Y$ and $g:B\to Y$,
consider whether or not there exists the homotopy commutative diagram
$$(5)\ \ \ \  
\xymatrix{E
\ar@{.>}[dr]^{F}& X\ar[l]_{i}\ar[d]^{f}\\
B \ar[r]_{g}\ar@{.>}[u]^s
 & Y\\
}$$
in which 
$\eta : X\overset{i}{\to} E\overset{p}{\to} B$
is  a fibration
with the section $s$.

Recall that a fibration $Z\overset{i}{\to} X\overset{p}{\to} Y$
is said to be {\it tncz} (totally non-cohomologous to zero) 
if 
$H^*(X)\cong H^*(Z)\otimes H^*(Y)$
as $H^*(Y)$-modules
by $p^*:H^*(Y)\to H^*(X)$. \\

\noindent
{\bf Definition B.} 
For a  map $f:X\to Y$, the $n$th {\it tncz
 subgroup} of $f$  is given as 
$${\mathcal T}_n(Y,X;f):=\{ a\in \pi_n(Y)\ |\ there \ are \ 
a \ tncz \ fibration  \ 
\eta  :X\overset{i}{\to} E\overset{p}{\to} S^n \ with$$
$$ a \ section \ s \ and \ a \ map \ F:E\to Y\
 such \ that\ F\circ s\simeq a,\  f\simeq  F\circ i 
\}.$$
\vspace{0.1cm}

\noindent
{\bf Definition C.} 
For a  map $f:X\to Y$, the $n$th {\it sectional 
 subgroup} of $f$  is given as 
$${\mathcal S}_n(Y,X;f):=\{ a\in \pi_n(Y)\ |\ there \ are \ 
a \ fibration  \ 
\eta  :X\overset{i}{\to} E\overset{p}{\to} S^n \ with$$
$$ a \ section \ s \ and \ a \ map \ F:E\to Y\
 such \ that\ F\circ s\simeq a,\  f\simeq  F\circ i 
\}.$$
\vspace{0.2cm}

\noindent
{\bf Remark 2.}  Note that ${\mathcal S}_n(Y,X;f)$ is a group.
For $a,b\in {\mathcal S}_n(Y,X;f)$,
there is a homotopy commutative diagram
$$
\xymatrix{\ar @{} [dr] |{\mbox{pull-back}}
E_{a+b}\ar[r]^{h}\ar[d]&E_{a}\cup_X E_b\ar[d]
\ar@{.>}[dr]^{G}& X\ar[l]_{\ \ i}\ar[d]^{f}\\
S^n\ar[r]_q&S^n\vee S^n \ar[r]_{(a|b)}
 & Y, \\
}$$
where $q$ is the pinching comultiplication and 
the dotted arrow $G$ is given  by the universality of push-out
from $F_{a,f}:E_a\to Y$ and $F_{b,f}:E_b\to Y$.
A section $s: S^n\to E_{a+b}$ is given by
$s(x):=(x,s'q(x))$ for a section $s':S^n\vee S^n\to E_{a}\cup_X E_b$
with $G\circ s'\simeq (a|b)$.
Also $i_{a+b}:X\to E_{a+b}$ is given by the universality of pull-back from $i:X\to E_a\cup_XE_b$ 
and $*:X\to S^n$.
Thus we have a homotopy commutative diagram
$$  
\xymatrix{E_{a+b}
\ar@{.>}[dr]^{G\circ h}& X\ar[l]_{i_{a+b}}\ar[d]^{f}\\
S^n \ar[r]_{a+b}\ar[u]^s
 & Y. \\
}$$
That is $a+b:=(a|b)\circ q\in {\mathcal S}_n(Y,X;f)$.
When $\eta_a:X\to E_a\to S^n$ and  $\eta_b:X\to E_b\to S^n$
are trivial (tncz),
the fibration $X\to E_{a+b}\to S^n$
 is trivial (tncz).
 Thus 
${\mathcal G}_n(Y,X;f)$
(${\mathcal T}_n(Y,X;f)$) is a group too.\\

We have the sequence of inclusions of groups:
$$G_n(f)
\subset {\mathcal G}_n(f)
\subset {\mathcal T}_n(f)\subset {\mathcal S}_n(f)\subset \pi_n(Y)$$
for a map $f:X\to Y$.

We consider them under some conditions on $X_{\Q}$
for a map $f:X\to Y$.
For the KS-extension of the fibration $\eta
:X\to E\to S^n$ in $(5)$,
we see $D'-d_X=0$  if $\pi_{\geq n}(X)\otimes \Q=0$.
Thus 

\begin{prop}
If  $\pi_{\geq n}(X)\otimes \Q=0$,
$ {\mathcal S}_n(Y_{\Q},X_{\Q};f_{\Q})=G_n(Y_{\Q},X_{\Q};f_{\Q})$.
\end{prop}

For example,  ${ G}_6(f_{\Q})= {\mathcal G}_6(f_{\Q})
={\mathcal T}_6(f_{\Q})= {\mathcal S}_6(f_{\Q})=0$
for the map $f$ in Example \ref{c}.
We know that
a fibration $X\to E\to B$ is homotopically trivial
if the classifying map $B\to Baut_1X$ is homotopic to the constant map.
Therefore we have

\begin{prop}
If $\pi_n(Baut_1X)\otimes \Q=0$,
 ${\mathcal G}_n(Y_{\Q},X_{\Q};f_{\Q})={ \mathcal S}_n(Y_{\Q},X_{\Q};f_{\Q})$.
\end{prop}

 In Example \ref{1} and Example \ref{2} below, 
  $\pi_n(Baut_1X)\otimes \Q=H_{n-1}(Der M(X))\neq 0$.
  But it is known that
  any fibration with fiber $\C P^m$, the $m$ dimensional complex projective space,
  is rationally tncz.
  In general, we note 
  
 \begin{lem}
  For $n>1$,
any fibration with fiber $X$ over $S^n$
  is rationally tncz if and only if 
  the  map $\rho :H_{n-1}(Der M(X))\to Der_{n-1} H^*(X;\Q)$
with  $\rho  ([\sigma ])([w])=[\sigma (w) ]$ is zero.
Here $Der_{n-1} H^*(X;\Q)$ means the  derivations of the graded algebra
$H^*(X;\Q)$ decreasing the degree by 
$n-1$ (\cite[9.7.2]{FOT}).
\end{lem}
  

\noindent
{\it Proof.}
The KS-extension of a fibration $X\to E\to S^n$
is given by the differential
$Dv=dv+\sigma (v)x$
for some $[\sigma ]\in H_{n-1}(Der M(X))$
with  the differential $d$ of $M(X)$
and $v\in M(X)$.
Then an element
 $[w]$ of $H^*(X;\Q)$ is extend to an element $[w+w'x]$ of   $H^*(E;\Q)$
 if and only if 
 $dw'=\sigma (w)$.
 \hfill\qed\\

\begin{exmp}\label{1}
{\rm For the associated fibration $S^2\to \C P^3\overset{f}{\to} S^4$ of 
the Hopf fibration $S^3\to S^7\overset{f}{\to} S^4$,
the KS-model is given by
$$(\Lambda 
(w_4,w_7),d_Y)\to (\Lambda (w_4,w_7,v_2,v_3),D) \to (\Lambda (v_2,v_3),d)$$
where $d_Yw_4=0$, $d_Yw_7=w_4^2$, $Dv_2=dv_2=0$,
$dv_3=v_2^2$ and $Dv_3=v_2^2-w_4$. 
Also $M(\C P^3)\cong 
(\Lambda (v_2,w_7),d_X)$ with $d_Xv_2=0$, $d_Xw_7=v_2^4$ and then 
$M(f)(w_4)=v_2^2$,
$M(f)(w_7)=w_7$.
Then we have $$ {\mathcal T}_4(S^4_{\Q},\C P^3_{\Q};f_{\Q})=\pi_4(S^4)\otimes \Q=\Q.$$

Indeed, 
for $a\in \pi_4(S^4)\otimes \Q$
with $M(a)(w_4)=cx$ ($c\in \Q$) and  $M(a)(w_7)=0$, 
put $$D'v_2=0,\ \ \  D'w_7=v_2^4+2cv_2^2x$$
$$\mbox{and }\ \ \ F(w_4)=v_2^2+cx, \ \ \ F(w_7)=w_7$$
 in
$$
\xymatrix{
(\Lambda x/x^2\otimes  \Lambda (v_2,w_7),{D'}) 
\ar[r]^{}\ar[d]_{} & 
(\Lambda (v_2,w_7 ),d_X)  \\
(\Lambda x/x^2,0)  &  
( \Lambda (w_4,w_7),d_Y).\ar[u]_{M(f)}\ar[l]^{M(a)}\ar@{.>}[ul]^{F}\\
}$$
Thus 
$a\in {\mathcal T}_4(S^4_{\Q},\C P^3_{\Q};f_{\Q})$.
On the other hand,  $ {\mathcal G}_4(S^4_{\Q},\C P^3_{\Q};f_{\Q})=0$
since 
$(\Lambda x/x^2\otimes  \Lambda (v_2,w_7),{D'})$ can  not be isomorphic to 
$ (\Lambda x/x^2,0)\otimes  (\Lambda (v_2,w_7),d_X) $ 
over $(\Lambda x/x^2,0)$ for any $D'$.
}
\end{exmp}

\begin{exmp}\label{2}{\rm 
For the map $f:\C P^2\to S^4$ collapsing the 2-cell,
 $M(f):M(S^4)= ( \Lambda (w_4,w_7),d_Y)
\to ( \Lambda (v_2,v_5),d_X)=M(\C P^2)$
is given by 
$M(f)(w_4)=v_2^2$ and $M(f)(w_7)=v_2v_5$.
Then we have
 $${\mathcal T}_4(S^4_{\Q},\C P^2_{\Q};f_{\Q})= \pi_4(S^4)\otimes \Q=\Q.$$
 
Indeed,
 for  $a\in \pi_4(S^4)\otimes \Q$
with $M(a)(w_4)=cx$  ($c\in \Q$) and  $M(a)(w_7)=0$,  
put
$$D'(v_2)=0,\ \ D'(v_5)=v_2^3+2cv_2x$$
$$\mbox{and }\ \ F(w_4)=v_2^2+cx,\ \ F(w_7)=v_2v_5$$
in
$$
\xymatrix{
(\Lambda x/x^2\otimes  \Lambda (v_2,v_5),{D'}) 
\ar[r]^{}\ar[d]_{} & 
(\Lambda (v_2,v_5 ),d_X)  \\
(\Lambda x/x^2,0)  &  
( \Lambda (w_4,w_7),d_Y).\ar[u]_{M(f)}\ar[l]^{M(a)}\ar@{.>}[ul]^{F}\\
}$$
Thus 
$a\in {\mathcal T}_4(S^4_{\Q},\C P^2_{\Q};f_{\Q})$.
On the other hand, we have
 ${\mathcal G}_4(S^4_{\Q},\C P^2_{\Q};f_{\Q})=0$ since 
$(\Lambda x/x^2\otimes  \Lambda (v_2,v_5),{D'})$ can  not be isomorphic to 
$ (\Lambda x/x^2,0)\otimes  (\Lambda (v_2,v_5),d_X) $ over $(\Lambda x/x^2,0)$
for any $D'$.
 
}
\end{exmp}

\begin{exmp}{\rm
 Put $\xi :\Omega X\overset{i}{\to} LX\overset{p}{\to} X$
 the fibration of free loops,
 in which $\Omega X$ is the loop space and $LX=map(S^1,X)$ is the free loop space of 
a simply connected space $X$.
It has the  section $s:X\to LX$ with $s(z)$ the constant loop at a  point $z$ in $X$. 
Consider the case that $X=S^{2}$.
Then ${\mathcal S}_2(LS^2,\Omega S^2;i)\ni s$ 
since we can choose  $F=id_{LS^2}$ as
 $$
\xymatrix{LS^2
\ar@{.>}[dr]^{=}& \Omega S^2\ar[l]_{i}\ar[d]^{i}\\
S^2 \ar[r]_{s}\ar[u]^s
 & LS^2.\\
}$$
Thus we have 
${\mathcal S}_2(LS^2,\Omega S^2;i)\neq 0$.
Especially,
we see ${\mathcal S}_2(LS^2_{\Q},\Omega S^2_{\Q};i_{\Q})\neq 0$
since $s$ is the torsion free generator of $\pi_2(LS^2)$.

But ${\mathcal T}_2(LS^2_{\Q},\Omega S^2_{\Q};i_{\Q})=0$.
Indeed,
the KS-model of $\xi$ is given by
$$(\Lambda (x,y),d_Y)\to 
(\Lambda (x,y,\overline{x},\overline{y}),D)\to(\Lambda (\overline{x},\overline{y}),0)$$
where $M(S^2)=(\Lambda (x,y),d_Y)$ with @$|x|=2$, $|y|=3$,
$d_Yx=0$, $d_Yy=x^2$, $|\overline{x}|=1$, $|\overline{y}|=2$,
$D(\overline{x})=0$, 
and $D(\overline{y})=2x\overline{x}$ \cite{VS}.
For the KS-model of a fibration
 $\eta  :\Omega S^2\overset{i}{\to} E\overset{p}{\to} S^2$
 with a section is given as
 $$(\Lambda (x,y),d_Y)\to 
(\Lambda (x,y,\overline{x},\overline{y}),D')\to(\Lambda (\overline{x},\overline{y}),0)$$
where  $D'(\overline{x})=0$ and 
$D'(\overline{y})=cx\overline{x}$ for some $c\in \Q$ by the degree arguments.
If $c=0$, there does not exist a map 
$F:(\Lambda (x,y,\overline{x},\overline{y}),D)\to (\Lambda (x,y,\overline{x},\overline{y}),D')=
(\Lambda (x,y,\overline{x},\overline{y}),d_Y)$ that we want.
If $c\neq 0$, it is not rationally  tncz since
$H^*(\Omega S^2;\Q )= \Lambda (\overline{x},\overline{y})$ and
$$H^*(E;\Q)\cong \Q [x,\{ u_i\}_{i>0} ]\otimes \Lambda (\overline{x})/(x^2,x\overline{x}, 
\{ xu_i\}, \{ \overline{x} u_i\},\{u_iu_j\} ),$$
where $u_i=[\overline{x} {\overline{y}}^i]$. 
}
\end{exmp}


\begin{exmp}{\rm
Put $\xi :S^1\overset{j}{\to} K\overset{f}{\to} S^1$
the fiber bundle with total space  a Klein bottle $K$.
Then
$G_1(K,S^1,j)={\mathcal G}_1(K,S^1,j)={\mathcal T}_1(K,S^1,j)=\Z$
and $ {\mathcal S}_1(K,S^1,j)=\pi_1 (K)$.
}
\end{exmp}


}
\end{document}